\documentclass{amsart}
\usepackage{amssymb,amsmath,latexsym}
\usepackage{times}
\usepackage[dvips]{graphicx}

\newtheorem{thm}{Theorem}[section]

\newtheorem{cor}[thm]{Corollary}
\newtheorem{prop}[thm]{Proposition}
\newtheorem{lemma}[thm]{Lemma}

\numberwithin{equation}{section}

\newcommand{\ssp}{\hspace{1pt}}
\newcommand{\del}{\backslash}
\newcommand{\cl}{\hbox{\rm cl}}
\newcommand{\mcZ}{\mathcal{Z}}
\newcommand{\mcF}{\mathcal{F}}
\newcommand{\mcY}{\mathcal{Y}}
\newcommand{\mcA}{\mathcal{A}}
\newcommand{\mcP}{\mathcal{P}}

\newcommand{\frp}{\mathbin{\Box}}

\title[Characterizations of Fundamental Transversal Matroids]
{Characterizations of Transversal and \\ Fundamental Transversal
  Matroids} \date{\today}
\author[J.~Bonin]{Joseph E.~Bonin}
\address{Department of Mathematics\\ The George Washington University\\
  Washington, D.C. 20052, USA} \email{jbonin@gwu.edu}
\author[J.~Kung]{Joseph P.\,S.~Kung}
\address{Department of Mathematics\\  University of North Texas\\
  Denton, TX 76203, USA} \email{kung@unt.edu}
\author[A.~de Mier]{Anna de Mier} \address{Departament de Matem\`atica
  Aplicada II\\ Universitat Polit\`ecnica de Catalunya\\ Jordi Girona
  1--3, 08034\\ Barcelona, Spain} \email{anna.de.mier@upc.edu}

\begin{document}

\begin{abstract}
  A result of Mason, as refined by Ingleton, characterizes transversal
  matroids as the matroids that satisfy a set of inequalities that
  relate the ranks of intersections and unions of nonempty sets of
  cyclic flats.  We prove counterparts, for fundamental transversal
  matroids, of this and other characterizations of transversal
  matroids.  In particular, we show that fundamental transversal
  matroids are precisely the matroids that yield equality in Mason's
  inequalities and we deduce a characterization of fundamental
  transversal matroids due to Brylawski from this simpler
  characterization.
\end{abstract}

\maketitle

\markboth{J.~Bonin, J.{\ssp}P.{\ssp}S.~Kung, A.~de Mier,
  \emph{Characterizations of Fundamental Transversal Matroids}}
{J.~Bonin, J.{\ssp}P.{\ssp}S.~Kung, A.~de Mier,
  \emph{Characterizations of Fundamental Transversal Matroids}}

\section{Introduction}

Transversal matroids can be thought of in several ways.  By
definition, a matroid is transversal if its independent sets are the
partial transversals of some set system.  A result of Brylawski gives
a geometric perspective: a matroid is transversal if and only if it
has an affine representation on a simplex in which each union of
circuits spans a face of the simplex.

Unions of circuits in a matroid are called \emph{cyclic sets}.  Thus,
a set $X$ in a matroid $M$ is cyclic if and only if the restriction
$M|X$ has no coloops.  Let $\mcZ(M)$ be the set of all cyclic flats of
$M$.  Under inclusion, $\mcZ(M)$ is a lattice: for $X,Y\in\mcZ(M)$,
their join in $\mcZ(M)$ is their join, $\cl(X\cup Y)$, in the lattice
of flats; their meet in $\mcZ(M)$ is the union of the circuits in
$X\cap Y$.  The following characterization of transversal matroids was
first formulated by Mason~\cite{mason} using sets of cyclic sets; the
observation that his result easily implies its streamlined counterpart
for sets of cyclic flats was made by Ingleton~\cite{ing}.  Theorem
\ref{thm:mi} has proven useful in several recent
papers~\cite{int,trl,cyc}.  For a family $\mcF$ of sets we shorten
$\cap_{X\in \mcF} X$ to $\cap \mcF$ and $\cup_{X\in \mcF} X$ to $\cup
\mcF$.

\begin{thm}\label{thm:mi}
  A matroid is transversal if and only if for all nonempty sets $\mcF$
  of cyclic flats,
  \begin{equation}\label{eq:mi}
    r(\cap \mcF) \leq
    \sum_{\mcF'\subseteq \mcF} (-1)^{|\mcF'|+1}r(\cup\mcF').
  \end{equation}
\end{thm}

It is natural to ask: which matroids satisfy the corresponding set of
equalities?  We show that $M$ satisfies these equalities if and only
if it is a fundamental transversal matroid, that is, $M$ is
transversal and it has an affine representation on a simplex (as
above) in which each vertex of the simplex has at least one matroid
element placed at it. The main part of this paper,
Section~\ref{sec:charftm}, provides four characterizations of these
matroids.

We recall the relevant preliminary material in Section \ref{sec:back}.
Theorems \ref{thm:ftm} and \ref{thm:ftmbeta} give new
characterizations of fundamental transversal matroids; from the
former, two other new characterizations (Theorem \ref{thm:ingthm3} and
Corollary \ref{cor:tisft}) follow easily.  The proofs of Theorems
\ref{thm:ftm} and \ref{thm:ftmbeta} use a number of ideas from a
unified approach to Theorem \ref{thm:mi} and a second characterization
of transversal matroids (the dual of another result of Mason,
from~\cite{masonalpha}); we present this material in
Section~\ref{sec:chartm} and deduce another of Mason's results from
it.  We conclude the paper with a section of observations and
applications; in particular, we show that Brylawski's characterization
of fundamental transversal matroids~\cite[Proposition 4.2]{aff}
follows easily from the dual of Theorem \ref{thm:ftm}.

As is common, we assume that matroids have finite ground sets.
However, no proofs use finiteness until we apply duality in
Theorem~\ref{thm:dualftm}, so, as we spell out in
Section~\ref{sec:app}, most of our results apply to matroids of finite
rank on infinite sets.

We assume basic knowledge of matroid theory; see~\cite{ox,welsh}.  Our
notation follows~\cite{ox}.  A good reference for transversal matroids
is~\cite{bru}.

We use $[r]$ to denote the set $\{1,2,\ldots,r\}$.

\section{Background}\label{sec:back}

Recall that a \emph{set system} $\mcA$ on a set $S$ is a multiset of
subsets of $S$.  It is convenient to write $\mcA$ as
$(A_1,A_2,\ldots,A_r)$ with the understanding that
$(A_{\sigma(1)},A_{\sigma(2)},\ldots,A_{\sigma(r)})$, where $\sigma$
is any permutation of $[r]$, is the same set system.  A \emph{partial
  transversal} of $\mcA$ is a subset $I$ of $S$ for which there is an
injection $\phi:I\rightarrow [r]$ with $x\in A_{\phi(x)}$ for all
$x\in I$.  \emph{Transversals} of $\mcA$ are partial transversals of
size $r$.  Edmonds and Fulkerson~\cite{ef} showed that the partial
transversals of a set system $\mcA$ on $S$ are the independent sets of
a matroid on $S$; we say that $\mcA$ is a \emph{presentation} of this
\emph{transversal matroid} $M[\mcA]$.

Of the following well-known results, all of which enter into our work,
Corollary \ref{cor:bmot} plays the most prominent role.  The proofs of
some of these results can be found in~\cite{bru}; the proofs of the
others are easy exercises.

\begin{lemma}\label{lem:r}
  Any transversal matroid $M$ has a presentation with $r(M)$ sets.  If
  $M$ has no coloops, then each presentation of $M$ has exactly $r(M)$
  nonempty sets.
\end{lemma}

\begin{lemma}
  If $M$ is a transversal matroid, then so is $M|X$ for each
  $X\subseteq E(M)$.  If $(A_1,\ldots,A_r)$ is a presentation $M$,
  then $(A_1\cap X,\ldots,A_r\cap X)$ is a presentation of $M|X$.
\end{lemma}

\begin{cor}\label{cor:bmot}
  If $(A_1,A_2,\ldots,A_r)$ is a presentation of $M$, then for each
  $F\in\mcZ(M)$, there are exactly $r(F)$ integers $i$ with $F\cap
  A_i\ne \emptyset$.
\end{cor}

\begin{lemma}
  For each $A_i\in \mcA$, its complement $A_i^c=E(M)-A_i$ is a flat of
  $M[\mcA]$.
\end{lemma}

\begin{lemma}\label{lem:extcy}
  If $(A_1,A_2,\ldots,A_r)$ is a presentation of $M$ and if $x$ is a
  coloop of $M\del A_i$, then $(A_1,A_2,\ldots,A_{i-1},A_i\cup
  x,A_{i+1},\ldots,A_r)$ is also a presentation of $M$.
\end{lemma}

\begin{cor}\label{lem:cycok}
  For any presentation $(A_1,A_2,\ldots,A_r)$ of a transversal matroid
  $M$, there is a presentation $(A'_1,A'_2,\ldots,A'_r)$ of $M$ with
  $A_i\subseteq A'_i$ and $A'^c_i\in\mcZ(M)$ for $i\in[r]$.
\end{cor}

A presentation $(A_1,A_2,\ldots,A_r)$ of $M$ is \emph{maximal} if,
whenever $(A'_1,A'_2,\ldots,A'_r)$ is a presentation of $M$ with
$A_i\subseteq A'_i$ for $i\in [r]$, then $A_i= A'_i$ for $i\in [r]$.
It is well known that each transversal matroid of rank $r$ has a
unique maximal presentation with $r$ sets.

A \emph{fundamental transversal matroid} is a transversal matroid that
has a presentation $(A_1,A_2,\ldots,A_r)$ for which no difference
$A_i-\bigcup_{j\in[r]-i}A_j$, for $i\in [r]$, is empty.  Clearly any
transversal matroid can be extended to a fundamental transversal
matroid: whenever a set in a given presentation is contained in the
union of the others, adjoin a new element to that set and to the
ground set, but to no other set in the presentation.

In the next paragraph, we describe how, given a presentation
$\mcA=(A_1,A_2,\ldots,A_r)$ of a transversal matroid $M$, we get an
affine representation of $M$ on a simplex.  (Our abstract description
is based on~\cite{aff}, which also discusses coordinates.)  Recall
that a simplex $\Delta$ in $\mathbb{R}^{r-1}$ is the convex hull of
$r$ affinely independent vectors, $v_1,v_2,\ldots,v_r$.  We will use
the following notation.  Given $\mcA$ and $\Delta$, for $x\in E(M)$,
let $\Delta_{\mcA}(x) = \{v_k\,:\,x\in A_k\}$; also, for $X\subseteq
E(M)$, let $\Delta_{\mcA}(X)=\cup_{x\in X}\Delta_{\mcA}(x)$.  (We omit
the subscript $\mcA$ when only one presentation is under discussion.)
The sets $\Delta(x)$ and $\Delta(X)$ span, and so can be identified
with, faces of $\Delta$.  Note that if $F\in\mcZ(M)$, then
$|\Delta(F)|=r(F)$ by Corollary \ref{cor:bmot}, so
$F=\{x\,:\,\Delta(x)\subseteq\Delta(F)\}$.

Given a presentation $\mcA$ of $M$, to get the corresponding affine
representation, first extend $M$ to a fundamental transversal matroid
$M'$ by extending $\mcA$ to a presentation $\mcA'$ of $M'$, as above.
We get an affine representation of $M'$ by, for each $x\in E(M')$,
placing $x$ as freely as possible (relative to all other elements) in
the face $\Delta_{\mcA'}(x)$ of $\Delta$.  Thus, a cyclic flat $F$ of
$M'$ of rank $i$ is the set of elements in some face of $\Delta$ with
$i$ vertices.  The affine representation of $M$ is obtained by
restricting that of $M'$ to $E(M)$.  Note that, by construction, such
an affine representation of $M$ can be extended to an affine
representation of a fundamental transversal matroid by adding elements
at the vertices of $\Delta$.

\begin{figure}
\begin{center}
  \includegraphics[width = 3.25truein]{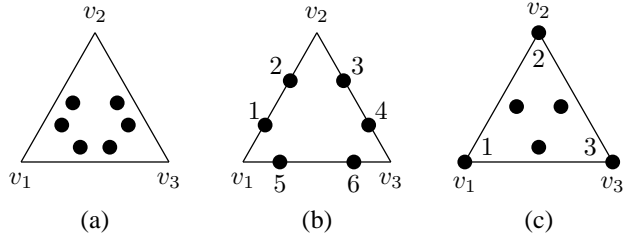}
\end{center}
\caption{Three representations of the uniform matroid $U_{3,6}$ on the
  simplex with vertices $v_1,v_2,v_3$.}\label{affrep}
\end{figure}

Such representations of the uniform matroid $U_{3,6}$ for the
presentations (a)~$([6],[6],[6])$,
(b)~$(\{1,2,5,6\},\{1,2,3,4\},\{3,4,5,6\})$, and
(c)~$(\{1,4,5,6\},\{2,4,5,6\},\{3,4,5,6\})$ are shown in
Figure~\ref{affrep}.  (Only elements $x$ with
$\Delta(x)\ne\{v_1,v_2,v_3\}$ are labelled in the figure.)

Note that the presentation can be recovered from the placement of the
elements.  The following result of Brylawski~\cite{aff} extends these
ideas.

\begin{thm}\label{thm:geomtr}
  A matroid $M$ is transversal if and only if it has an affine
  representation on a simplex $\Delta$ in which, for each $F\in
  \mcZ(M)$, the flat $F$ is the set of elements in some face of
  $\Delta$ with $r(F)$ vertices.
\end{thm}

With this result, we can give a second perspective on fundamental
transversal matroids.  A basis $B$ of a matroid $M$ is a
\emph{fundamental basis} if each $F\in\mcZ(M)$ is spanned by $B\cap
F$.  In any affine representation of a matroid $M$ with a fundamental
basis $B$, if the elements of $B$ are placed at the vertices of a
simplex $\Delta$, then a cyclic flat of rank $i$ is the set of
elements in some $i$-vertex face of $\Delta$.  It follows from Theorem
\ref{thm:geomtr} that a matroid is a fundamental transversal matroid
if and only if it has a fundamental basis.

We use the following terminology from ordered sets, applied to the
lattice $\mcZ(M)$ of cyclic flats.  An \emph{antichain} in $\mcZ(M)$
is a set $\mcF\subseteq \mcZ(M)$ such that no two sets in $\mcF$ are
related by inclusion.  A \emph{filter} in $\mcZ(M)$ is a set
$\mcF\subseteq\mcZ(M)$ such that if $A\in \mcF$ and $B\in \mcZ(M)$,
and if $A\subseteq B$, then $B\in \mcF$.  An \emph{ideal} in $\mcZ(M)$
is a set $\mcF\subseteq\mcZ(M)$ such that if $B\in \mcF$ and $A\in
\mcZ(M)$, and if $A\subseteq B$, then $A\in \mcF$.

\section{Characterizations of Transversal Matroids}\label{sec:chartm}

In the main part of this section, we connect Theorem \ref{thm:mi} with
another characterization of transversal matroids by giving a cycle of
implications that proves both.  While parts of the argument have
entered into proofs of related results, the link between these results
seems not to have been exploited before.  In Section~\ref{sec:charftm}
we use substantial parts of the material developed here.  We end this
section by showing how another characterization of transversal
matroids follows easily from Theorem \ref{thm:mi}.

To motivate the second characterization (part (3) of Theorem
\ref{thm:betami}), we describe how to prove that a matroid $M$ that
satisfies the condition in Theorem \ref{thm:mi} is transversal.  We
want to construct a presentation of $M$.  By Corollary
\ref{lem:cycok}, $M$ should have a presentation $\mcA$ in which the
complement of each set $A_i$ is in $\mcZ(M)$.  Thus, we must
determine, for each $F\in\mcZ(M)$, the multiplicity of $F^c$ in
$\mcA$.  We will define a function $\beta$ on all subsets of $E(M)$ so
that for each $F\in\mcZ(M)$, the multiplicity of $F^c$ in $\mcA$ will
be $\beta(F)$.  In particular, the sum of $\beta(F)$ over all $F\in
\mcZ(M)$, i.e., $|\mcA|$, should be $r(M)$.  By Corollary
\ref{cor:bmot}, for each $F\in \mcZ(M)$ we must have
\begin{equation}\label{eq:bc}
\sum_{Y\in\mcZ(M)\,:F\cap Y^c\ne\emptyset}\beta(Y) = r(F),
\end{equation}
or, equivalently,
\begin{equation}\label{eq:betacy}
  \sum_{Y\in\mcZ(M)\,:F\subseteq Y}\beta(Y) = r(M)-r(F).
\end{equation}
With this motivation, we define $\beta$ recursively on all subsets $X$
of $E(M)$ by
\begin{equation}\label{eq:beta}
  \beta(X) = r(M)-r(X)- \sum_{Y\in \mcZ(M)\,:\, X\subset Y}\beta(Y). 
\end{equation} 
By the definition of $\beta$, equation (\ref{eq:betacy}) holds
whenever $F$ spans a cyclic flat of $M$.  Applying that equation to
the cyclic flat $\cl(\emptyset)$ gives
\begin{equation}\label{eq:betasum}
  \sum_{Y\in\mcZ(M)} \beta(Y) = r(M).
\end{equation} 
Thus, equation (\ref{eq:bc}) follows for $F \in \mcZ(M)$.

(The function $\beta$ is dual to the function $\alpha$ that was
introduced in~\cite{masonalpha} and studied further
in~\cite{ip,alpha}; see the comments in the first part of
Section~\ref{sec:app}.  The definition of the function $\tau$
in~\cite{bru} is similar to that of $\beta$, although values of $\tau$
that would otherwise be negative are set to zero; with the recursive
nature of the definition, this can change the values on more sets than
just those on which $\beta$ is negative.  It follows from
Theorem~\ref{thm:betami} that $\beta$ and $\tau$ agree precisely on
transversal matroids.)

The next lemma plays several roles.

\begin{lemma}\label{lem:keysum}
  If $\mcF$ is a nonempty filter in $\mcZ(M)$, then
  \begin{equation}\label{betasum}
    \sum_{Y\in\mcF}\beta(Y) = r(M) - \sum_{\mcF'\,:\, \mcF' \subseteq \mcF}
    (-1)^{|\mcF'|+1}\, r(\cup\mcF').
  \end{equation}
  Also, if $\mcF_0$ is any subset of $\mcF$ that contains every
  minimal set in $\mcF$, then the sum on the right can be taken just
  over all nonempty subsets $\mcF'$ of $\mcF_0$.
\end{lemma}

\begin{proof}
  For each $Y\in \mcF$, the set $\mcY=\{F\in\mcF\,:\,F \subseteq Y\}$
  is nonempty, so
  \begin{equation*}
    \sum_{\mcF'\subseteq \mcY\,:\,\mcF' \ne\emptyset} (-1)^{|\mcF'|+1}
    = 1.
  \end{equation*}  
  From this sum and equation (\ref{eq:betacy}), we have
  \begin{align*}\label{ekey}
    \sum_{Y\in\mcF}\beta(Y) &\, = \sum_{Y\in\mcF}\beta(Y)
    \sum\limits_{\substack {\mcF'\,:\, \mcF' \subseteq
        \{F\in\mcF\,:\,F \subseteq Y\}, \\
        \mcF'\ne\emptyset}} (-1)^{|\mcF'|+1} \\
    &\, = \sum\limits_{\substack {\mcF'\,:\, \mcF' \subseteq \mcF, \\
        \mcF'\ne\emptyset}}
    (-1)^{|\mcF'|+1} \sum\limits_{Y\in\mcF\,:\, \cup\mcF'\subseteq Y} \beta (Y) \\
    &\, = \sum\limits_{\substack {\mcF'\,:\, \mcF' \subseteq \mcF, \\
        \mcF'\ne\emptyset}} (-1)^{|\mcF'|+1}\,
    \bigl(r(M)-r(\cup\mcF')\bigr).
  \end{align*}
  Simplification yields equation (\ref{betasum}).  To prove the second
  assertion, note that for $X, Y \in \mcF$ with $X\subset Y$, the
  terms $(-1)^{|\mcF'|+1}\, r(\cup\mcF')$ with $Y\in \mcF'$ cancel via
  the involution that adjoins $X$ to, or omits $X$ from, $\mcF'$.
\end{proof}

We now turn to the first two characterizations of transversal
matroids.  The last part of the proof of Theorem \ref{thm:betami} uses
Hall's theorem: a set system $\mcA$ with $r$ sets has a transversal if
and only if, for each $i\in[r]$, each union of $i$ sets in $\mcA$ has
at least $i$ elements.

\begin{thm}\label{thm:betami}
  For a matroid $M$, the following statements are equivalent:
  \begin{enumerate}
  \item $M$ is transversal,
  \item for every nonempty subset (equivalently, filter; equivalently,
    antichain) $\mcF$ of $\mcZ(M)$,
    \begin{equation}
      r(\cap \mcF) \leq \sum_{\mcF'\subseteq \mcF}
      (-1)^{|\mcF'|+1}r(\cup \mcF'),\label{mi2} 
    \end{equation}
  \item $\beta(X)\geq 0$ for all $X\subseteq E(M)$.
  \end{enumerate}
\end{thm}

\begin{proof}
  The three formulations of statement (2) are equivalent since, for
  $X, Y \in \mcF$ with $X\subset Y$, using $\mcF - \{Y\}$ in place of
  $\mcF$ preserves the right side of the inequality by the argument at
  the end of the proof of Lemma \ref{lem:keysum}; also, the left side
  is clearly the same.

  To show that statement (1) implies statement (2), extend $M$ to a
  fundamental transversal matroid $M_1$.  Let $r_1$ and $\cl_1$ be its
  rank function and closure operator.  For $\mcF\subseteq \mcZ(M)$,
  setting $\mcF_1 =\{\cl_1(F)\,:\,F\in \mcF\}$ gives
  $\mcF_1\subseteq\mcZ(M_1)$ as well as $r(\cup \mcF) = r_1(\cup
  \mcF_1)$ and $r(\cap \mcF) \leq r_1(\cap \mcF_1)$, so statement (2)
  will follow by showing that for $M_1$, equality holds in inequality
  (\ref{mi2}).

  Let $B$ be a fundamental basis of $M_1$ and let $\mcF\subseteq
  \mcZ(M_1)$ be nonempty.  We claim that
  \begin{equation}\label{eq:capcup}
    r_1(\cup \mcF)= \bigl|B\cap (\cup \mcF)\bigl| \qquad \text{and} \qquad
    r_1(\cap \mcF) = \bigl|B\cap (\cap \mcF)\bigr|.
  \end{equation}
  The first equality holds since $B\cap (\cup \mcF)$ is independent
  and each $F \in \mcF$ is spanned by $B\cap F$.  For the second
  equality, we have $r_1(\cap \mcF) \geq \bigl|B\cap (\cap
  \mcF)\bigr|$ since $B$ is independent.  To show that $B\cap (\cap
  \mcF)$ spans $\cap \mcF$, consider $x \in (\cap \mcF)-B$.  Since $x$
  is not in the basis $B$, the set $B\cup x$ contains a unique
  circuit, say $C$.  Clearly, $x\in C$.  Similarly using the basis
  $B\cap F$ of $F$, for $F\in \mcF$, and the uniqueness of $C$ gives
  $C-x\subseteq B\cap F$; thus, $C-x\subseteq B\cap (\cap \mcF)$, so,
  as needed, $B\cap (\cap \mcF)$ spans $\cap \mcF$.

  For $M_1$, upon using equations (\ref{eq:capcup}) to rewrite both
  sides of inequality (\ref{mi2}), it is easy to see that equality
  follows from inclusion-exclusion.

  We now show that statement (2) implies statement (3).  For
  $X\subseteq E(M)$, let $\mcF(X)$ be $ \{Y\in \mcZ(M)\,:\,X\subset
  Y\}$.  By equation (\ref{eq:beta}), proving $\beta(X)\geq 0$ is the
  same as proving
  \begin{equation*}
    \sum_{Y\in \mcF(X)}\beta(Y)\leq r(M)-r(X).
  \end{equation*}  
  This inequality is clear if $\mcF(X) = \emptyset$; otherwise, it
  follows from Lemma \ref{lem:keysum}, statement (2), and the obvious
  inequality $r(X)\leq r(\cap \mcF(X))$.
    
  Lastly, to show that statement (3) implies statement (1), we show
  that $M=M[\mcA]$ where $\mcA=(F^c_1,F^c_2, \ldots,F^c_r)$ is the
  multiset that consists of $\beta(F)$ occurrences of $F^c$ for each
  cyclic flat $F$ of $M$.  By equation (\ref{eq:betasum}), we have
  $r=r(M)$.

  To show that each dependent set $X$ of $M$ is dependent in
  $M[\mcA]$, it suffices to show this when $X$ is a circuit of $M$.
  In this case, $\cl_M(X)$ is a cyclic flat of $M$, so, by equation
  (\ref{eq:bc}) and the definition of $\mcA$, it has nonempty
  intersection with exactly $r_M(X)$ sets of $\mcA$, counting
  multiplicity.  Thus, $X$ is dependent in $M[\mcA]$ since, with
  $r_M(X)<|X|$, it cannot be a partial transversal of $\mcA$.

  To show that each independent set of $M$ is independent in
  $M[\mcA]$, it suffices to show this for each basis $B$.  For this,
  we use Hall's theorem to show that $(F^c_1\cap B,\ldots,F^c_r\cap
  B)$ has a transversal (which necessarily is $B$).  Let $X =
  \bigcup_{j\in J} (F^c_j\cap B)$ with $J\subseteq [r]$.  We must show
  $|X|\geq |J|$.  Now $B-X\subseteq\bigcap_{j\in J} F_j$, so
  \begin{align*}
    |J|\leq &\,\, \sum_{Y\in \mcZ(M)\,:\,B-X\subseteq Y}\beta(Y)\\
    \leq &\,\, r(M) -r(B-X)\\
    = &\,\, |B| - |B-X|\\
    \leq &\,\, |X|.
  \end{align*}
  (Note that in these inequalities, statement (3) is used twice in the
  first two lines.)
\end{proof}

It follows from equation (\ref{eq:betacy}) that the definition of
$\beta$ on cyclic flats is forced by wanting a presentation in which
the complement of each set is a cyclic flat.  Maximal presentations
have this property by Corollary \ref{lem:cycok}, so we have the next
result.

\begin{cor}\label{cor:uniquemax}
  The maximal presentation $\mcA$ of $M$ is unique; it consists of the
  sets $F^c$ with $F\in\mcZ(M)$, where $F^c$ has multiplicity
  $\beta(F)$ in $\mcA$.
\end{cor}

Like Theorem \ref{thm:mi}, the next result is a refinement, noted by
Ingleton~\cite{ing}, of a result of Mason~\cite{mason} that used
cyclic sets.  Mason used this result in his proof of Theorem
\ref{thm:mi}; we show that it follows easily from that result.  Let
$2^{[r]}$ be the lattice of subsets of $[r]$.

\begin{thm}\label{thm:ingthm2}
  A matroid $M$ of rank $r$ is transversal if and only if there is an
  injection $\phi:\mcZ(M)\rightarrow 2^{[r]}$ with
  \begin{enumerate}
  \item $|\phi(F)|=r(F)$ for all $F\in \mcZ(M)$,
  \item $\phi\bigl(\cl(F\cup G)\bigr)=\phi(F)\cup\phi(G)$ for all
    $F,G\in \mcZ(M)$, and
  \item $r(\cap\mcF)\leq |\cap\{\phi(F)\,:\,F\in\mcF\}|$ for every
    subset (equivalently, filter; equivalently, antichain) $\mcF$ of
    $\mcZ(M)$.
  \end{enumerate}
\end{thm}

\begin{proof}
  Assume $M = M[\mcA]$ with $\mcA = (A_1,A_2,\ldots,A_r)$.  For
  $F\in\mcZ(M)$, let
  \begin{equation}\label{phidef}
    \phi(F) = \{k\,:\,F\cap A_k\ne\emptyset\}.
  \end{equation}
  It is easy to see that $\phi$ is an injection and that properties
  (1)--(3) hold; in particular, the first is Corollary \ref{cor:bmot}.
  For the converse, assume $\phi:\mcZ(M)\rightarrow 2^{[r]}$ is an
  injection that satisfies properties (1)--(3).  For any nonempty
  subset $\mcF$ of $\mcZ(M)$, properties (1) and (2) allow us to
  recast the right side of inequality (\ref{mi2}) as the summation
  part of an inclusion-exclusion equation for the sets $\phi(F)$ with
  $F\in\mcF$; inequality (\ref{mi2}) follows from inclusion-exclusion
  and property (3), so $M$ is transversal by Theorem~\ref{thm:betami}.
\end{proof}

\section{Characterizations of Fundamental Transversal
  Matroids}\label{sec:charftm}

In this section, we treat counterparts, for fundamental transversal
matroids, of the results in the last section.  In contrast to
Theorem~\ref{thm:mi}, in the main result, Theorem~\ref{thm:ftm}, we
must work with cyclic flats since equality~(\ref{eq:mi3}) may fail for
sets $\mcF$ of cyclic sets.

\begin{thm}\label{thm:ftm}
  A matroid $M$ is a fundamental transversal matroid if and only if
  \begin{equation}\label{eq:mi3}
    r(\cap \mcF)=\sum_{\mcF'\subseteq \mcF}(-1)^{|\mcF'|+1}r(\cup \mcF')
  \end{equation}
  for all nonempty subsets (equivalently, antichains; equivalently,
  filters) $\mcF\subseteq\mcZ(M)$.
\end{thm}

In the proof of Theorem \ref{thm:betami}, we showed that equation
(\ref{eq:mi3}) holds for all fundamental transversal matroids; below
we prove the converse.  In the proof, we use the notation $\Delta(x)$
and $\Delta(X)$ that we defined in Section~\ref{sec:back}.  The
following well-known lemma is easy to prove.

\begin{lemma}\label{lem:delcir}
  If $C$ is a circuit of $M$, then $\Delta(C) = \Delta(C-x)$ for all
  $x\in C$.
\end{lemma}

By Theorem \ref{thm:betami}, if equation (\ref{eq:mi3}) always holds,
then $M$ is transversal.  In this setting, the next lemma identifies
$|\Delta(x)|$, for the maximal presentation, as the rank of a set.

\begin{lemma}\label{lem:int}
  Let $(A_1,A_2,\ldots,A_r)$ be the maximal presentation of a matroid
  $M$ for which equality (\ref{eq:mi3}) holds for all nonempty subsets
  of $\mcZ(M)$.  For each $x\in E(M)$, we have
  $|\Delta(x)|=r(\cap\mcF)$ where $\mcF=\{F\in\mcZ(M)\,:\,x\in F\}$.
\end{lemma}

\begin{proof}
  The set $\Delta(x)$ contains the vertices $v_k$ where $A_k = F^c$
  and $F\in\mcZ(M)-\mcF$.  By Lemma \ref{lem:keysum}, Corollary
  \ref{cor:uniquemax}, and equations (\ref{eq:betasum}) and
  (\ref{eq:mi3}), $|\Delta(x)|$ is, as stated,
  \begin{equation*}
    \sum_{F\in\mcZ(M)-\mcF}\beta(F) =r(\cap \mcF).\qedhere
  \end{equation*}
\end{proof}

The equality $|\Delta(x)|=r(\cap\mcF)$ may fail if equality
(\ref{eq:mi3}) fails.  For example, consider the rank-$4$ matroid on
$\{a,b,c,d,e,f,g\}$ in which $\{a,b,c,d\}$ and $\{d,e,f,g\}$ are the
only non-spanning circuits.  In the affine representation arising from
the maximal presentation, $d$ is placed freely on an edge of the
simplex even though the cyclic flats that contain it intersect in rank
one.

We now prove the main result.

\begin{proof}[Proof of Theorem \ref{thm:ftm}]
  Assume equation (\ref{eq:mi3}) holds for all nonempty sets of cyclic
  flats.  As noted above, $M$ is transversal.  Coloops can be placed
  at vertices of $\Delta$ and doing so reduces the problem to a
  smaller one, so we may assume that $M$ has no coloops.  Thus, $E(M)
  \in \mcZ(M)$.  The set $V$ of vertices of $\Delta$ has size $r(M)$,
  so $\Delta(E(M)) =V$.

  Let $\mcA$ be the maximal presentation of $M$.  As parts (a) and (c)
  of Figure~\ref{affrep} show, from the corresponding affine
  representation, it may be possible to get other affine
  representations of $M$ by moving some elements of $M$ to vertices of
  $\Delta$, where $x\in E(M)$ may be moved only to a vertex in
  $\Delta_{\mcA}(x)$.  Such affine representations correspond to
  presentations $\mcA'$ of $M$ in which, for each $x\in E(M)$, either
  $\Delta_{\mcA'}(x)=\Delta_{\mcA}(x)$ or $\Delta_{\mcA'}(x)=\{v_i\}$
  for some $v_i\in \Delta_{\mcA}(x)$.  Among all such affine
  representations, fix one with the minimum number of vertices of
  $\Delta$ at which no element of $E(M)$ is placed; let $\mcA'$ be the
  corresponding presentation.  To show that $M$ is fundamental, we
  show that if, in this affine representation, no element is placed at
  vertex $v_i$ of $\Delta$, then we get another affine representation
  of $M$ by moving some element there, which contradicts the
  minimality assumption.

  To show this, we will use the fundamental transversal matroid $M_1$
  that we obtain from the fixed affine representation of $M$
  (corresponding to $\mcA'$) by adding an element (which we call
  $v_j$) at each vertex $v_j$ of $\Delta$ at which there is no element
  of $M$.  Let $\mcP$ be the corresponding presentation of $M_1$.  Let
  $r_1$ and $\cl_1$ be its rank function and closure operator.  For
  $\mcF\subseteq \mcZ(M)$, let $\mcF_1 = \{\cl_1(F)\,:\,F\in\mcF\}$.
  Clearly $r(\cup\mcF) = r_1(\cup\mcF_1)$.  We claim that
  \begin{itemize}
  \item[(i)] $r(\cap \mcF)=r_1(\cap \mcF_1)$,
  \item[(ii)] $\Delta_{\mcA'}(\cap \mcF) = \Delta_{\mcP}(\cap
    \mcF_1)$, and
  \item[(iii)] $r(\cap\mcF)= |\Delta_{\mcA'}(\cap \mcF)|$.
  \end{itemize}
  To prove these properties, note that since $M_1$ is fundamental, we
  have
  \begin{equation*}
    \sum_{\mcF'\subseteq  \mcF_1}(-1)^{|\mcF'|+1}r_1(\cup \mcF')=
    r_1(\cap \mcF_1).
  \end{equation*}  
  Term by term, the sum agrees with its counterpart for $\mcF$ in $M$,
  so property (i) follows from equation (\ref{eq:mi3}).  Clearly,
  $r(\cap\mcF)\leq |\Delta_{\mcA'}(\cap \mcF)|$.  Also,
  $\Delta_{\mcA'}(\cap \mcF) \subseteq \Delta_{\mcP}(\cap \mcF_1)$
  since $\cap \mcF_1$ contains $\cap \mcF$.  Since $M_1$ is
  fundamental, equation (\ref{eq:capcup}) holds, from which we get
  $|\Delta_{\mcP}(\cap \mcF_1)|=r_1(\cap\mcF_1)$.  With these
  deductions, property (i) gives properties (ii) and (iii).

  Now assume that no element of $M$ has been placed at vertex $v_i$ of
  $\Delta$.  Let
  \begin{equation*}
    \mcF = \{F\in \mcZ(M)\,:\, v_i\in\Delta_{\mcA'}(F)\}.
  \end{equation*}  
  (By Corollary \ref{cor:bmot}, $\Delta_{\mcA}(F) = \Delta_{\mcA'}(F)$
  for all $F\in \mcZ(M)$.)  Now $E(M)\in\mcF$, so $\mcF\ne\emptyset$.
  Since $v_i \in \Delta_{\mcP}(\cap \mcF_1)$, property (ii) gives
  $v_i\in \Delta_{\mcA'}(\cap \mcF)$.  Fix $x \in \cap \mcF$ with
  $v_i\in\Delta_{\mcA'}(x)$.

  We claim that $\mcF = \{ F\in\mcZ(M)\,:\, x\in F\}$.  If $F\in
  \mcF$, then $\cap \mcF\subseteq F$, so $x\in F$.  Conversely, if
  $F\in\mcZ(M)$ and $x\in F$, then $v_i\in \Delta_{\mcA'}(F)$ since
  $\Delta_{\mcA'}(x)\subseteq \Delta_{\mcA'}(F)$.

  Now $v_i\in \Delta_{\mcA'}(x)$ but $x$ was not placed at $v_i$, so
  $\Delta_{\mcA'}(x) = \Delta_{\mcA}(x)$.  Since $x \in \cap \mcF$, we
  have $\Delta_{\mcA'}(x) \subseteq \Delta_{\mcA'}(\cap \mcF)$;
  property (iii), the previous paragraph, and Lemma \ref{lem:int} give
  equality, that is, $x$ is placed freely in the face
  $\Delta_{\mcA'}(\cap\mcF)$.  Let $M_2$ be the matroid that is
  obtained by moving $x$ to $v_i$, that is, $M_2 = M[\mcA'']$ where
  $\mcA''$ is formed from $\mcA'$ by removing $x$ from all sets except
  the one indexed by $i$.  We claim that $M$ and $M_2$ have the same
  circuits and so are the same matroid, thus proving our claim that
  some element can be moved to $v_i$.  Among all sets $C$ that are
  circuits of just one of $M$ and $M_2$ (if there are any), let $C$
  have minimum size.  Clearly, $x \in C$.

  We claim that $\Delta_{\mcA'}(C)=\Delta_{\mcA''}(C)$.  If $C$ is a
  circuit of $M$, then the claim follows from Lemma \ref{lem:delcir},
  the inclusion $\Delta_{\mcA''}(x)\subset\Delta_{\mcA'}(x)$, and the
  observation that $\Delta_{\mcA'}(y)=\Delta_{\mcA''}(y)$ for $y\in
  C-x$.  Assume $C$ is a circuit of $M_2$.  By Lemma \ref{lem:delcir},
  $v_i\in \Delta_{\mcA''}(y)$ for some $y$ in $C-x$.  Thus, $v_i\in
  \Delta_{\mcA'}(y)$, so all cyclic flats that contain $y$ are in
  $\mcF$ and so contain $x$; thus, all sets in the maximal
  presentation that contain $x$ also contain $y$, that is,
  $\Delta_{\mcA}(x)\subseteq \Delta_{\mcA}(y)$.  Since no element
  prior to $x$ was placed at $v_i$, we have $\Delta_{\mcA'}(y) =
  \Delta_{\mcA}(y)$; also, as noted above, $\Delta_{\mcA'}(x) =
  \Delta_{\mcA}(x)$, so $\Delta_{\mcA'}(x)\subseteq
  \Delta_{\mcA'}(y)$, from which the claim follows.
 
  Now $C$ is a circuit in one of $M$ and $M_2$, so, since
  $\Delta_{\mcA'}(C)=\Delta_{\mcA''}(C)$, we have
  \begin{equation*}  
    |\Delta_{\mcA'}(C)|=|\Delta_{\mcA''}(C)|<|C|.
  \end{equation*}  
  It follows that $C$ is dependent in both $M$ and $M_2$.  From this
  conclusion and the minimality assumed for $|C|$, it follows that $C$
  cannot be a circuit of just one of $M$ and $M_2$.  Thus, $M$ and
  $M_2$ have the same circuits and so are the same matroid, as we
  needed to show.
\end{proof}

The following result is immediate from Theorem \ref{thm:ftm} and
Lemma \ref{lem:keysum}.

\begin{thm}\label{thm:ftmbeta}
  A matroid $M$ is a fundamental transversal matroid if and only if
  \begin{equation}\label{eq:betasumfil}
    \sum_{Y\in\mcF}\beta(Y) = r(M) - r(\cap \mcF)
  \end{equation}
  for all filters $\mcF\subseteq \mcZ(M)$.
\end{thm}

The proof of the next result is similar to that of Theorem
\ref{thm:ingthm2} and uses Theorem \ref{thm:ftm}.

\begin{thm}\label{thm:ingthm3}
  A matroid $M$ of rank $r$ is a fundamental transversal matroid if
  and only if there is an injection $\phi:\mcZ(M)\rightarrow 2^{[r]}$
  with
  \begin{enumerate}
  \item $|\phi(F)|=r(F)$ for all $F\in \mcZ(M)$,
  \item $\phi\bigl(\cl(F\cup G)\bigr)=\phi(F)\cup\phi(G)$ for all
    $F,G\in \mcZ(M)$, and
  \item $r(\cap\mcF)=|\cap\{\phi(F)\,:\,F\in\mcF\}|$ for every subset
    (equivalently, filter; equivalently, antichain) $\mcF$ of
    $\mcZ(M)$.
  \end{enumerate}
\end{thm}

If the matroid $M$ is already known to be transversal and if a
presentation of $M$ is known, then we should define the function
$\phi$ in the last result as in equation (\ref{phidef}) or,
equivalently, $\phi(F) = \{k\,:\,v_k\in\Delta(F)\}$.  Properties (1)
and (2) then hold, so we have the next corollary.

\begin{cor}\label{cor:tisft}
  Let $\mcA$ be any presentation of a transversal matroid $M$.  The
  matroid $M$ is fundamental if and only if $r(\cap \mcF) =
  |\cap\{\Delta(F)\,:\,F\in\mcF\}|$ for every subset (equivalently,
  filter; equivalently, antichain) $\mcF$ of $\mcZ(M)$.
\end{cor}

\section{Observations and Applications}\label{sec:app}

We first consider the duals of the results above.  In particular,
Theorem \ref{thm:abdual} makes precise the remark before Lemma
\ref{lem:keysum}, that $\beta$ is dual to Mason's function $\alpha$;
this shows that the equivalence of statements (1) and (3) in Theorem
\ref{thm:betami} is the dual of Mason's result that $M$ is a
cotransversal matroid (a strict gammoid) if and only if $\alpha(X)\geq
0$ for all $X\subseteq E(M)$.

It is well known and easy to prove that
\begin{equation}\label{eq:zdual}
  \mcZ(M^*)=\{E(M)-F\,:\,F\in\mcZ(M)\},
\end{equation}
where $M^*$ is the dual of $M$.  With this result and the formula
\begin{equation}\label{eq:dual}
  r^*(X) = |X|-r(M)+r\bigl(E(M)-X\bigr)
\end{equation}
for the rank function $r^*$ of $M^*$, it is routine to show that a
matroid $M$ satisfies statement (2) in Theorem \ref{thm:betami} if and
only if for all sets (equivalently, ideals; equivalently, antichains)
$\mcF\subseteq \mcZ(M^*)$,
\begin{equation*}
  r^*(\cup \mcF) \leq \sum_{\mcF'\subseteq \mcF\,:\,\mcF'\ne\emptyset}
  (-1)^{|\mcF'|+1}r^*(\cap \mcF').
\end{equation*}
Thus, this condition characterizes cotransversal matroids $M^*$.

We now recall the function $\alpha$ that Mason introduced
in~\cite{masonalpha}, which is defined recursively as follows.  For
$X\subseteq E(M)$, set
\begin{equation}\label{alphadef}
  \alpha(X)=\eta(X)-\sum_{\text{flats } F\,:\,F\subset X} \alpha(F),
\end{equation}
where $\eta(X)$ is the nullity, $|X|-r(X)$, of $X$.  Thus for any flat
$X$ of $M$,
\begin{equation}\label{alphasum}
  \sum_{\text{flats } F\,:\,F\subseteq X} \alpha(F)=\eta(X).
\end{equation}

To prepare to link the functions $\alpha$ and $\beta$, we first show
that $\alpha(F)=0$ if $F$ is a noncyclic flat.  Induct on $|F|$.  The
base case holds vacuously.  Let $I$ be the set of coloops of $M|F$ and
set $F'=F-I$, so $F'\in\mcZ(M)$.  Since $\eta(F)=\eta(F')$, equation
(\ref{alphasum}) gives
\begin{equation*}
  \sum_{\text{flats } Y\,:\,Y\subseteq F} \alpha(Y)
  = \sum_{\text{flats } Y'\,:\,Y'\subseteq F'} \alpha(Y').
\end{equation*}  
Now $F$ and $F'$ contain precisely the same cyclic flats, so
$\alpha(F)$ is the only term in which the two sides of this equality
differ that is not yet known to be zero, so $\alpha(F)=0$.

It now follows that the sum in equation (\ref{alphadef}) can be over
just $F\in\mcZ(M)$ with $F\subset X$.  With induction, the next
theorem follows from this result and equations (\ref{eq:zdual}) and
(\ref{eq:dual}).

\begin{thm}\label{thm:abdual}
  For any matroid $M$, if $X\subseteq E(M)$, then
  $\alpha_M(X)=\beta_{M^*}\bigl(E(M)-X\bigr)$.
\end{thm}

As shown in~\cite{mlv}, the class of fundamental transversal matroids
is closed under duality.  (To see this, note that a basis $B$ of $M$
is fundamental if and only if $r(M)=r(F)+|B-F|$ for every
$F\in\mcZ(M)$; a routine rank calculation then shows that $B$ is a
fundamental basis of $M$ if and only if $E(M)-B$ is a fundamental
basis of $M^*$.)  Using this result and those above, it is easy to
deduce the following dual versions of Theorems \ref{thm:ftm} and
\ref{thm:ftmbeta}.  (Likewise, one can dualize Theorem
\ref{thm:ingthm3} and Corollary \ref{cor:tisft}.)

\begin{thm}\label{thm:dualftm}
  For a matroid $M$, the following statements are equivalent:
  \begin{enumerate}
  \item $M$ is a fundamental transversal matroid,
  \item for all subsets (equivalently, ideals; equivalently,
    antichains) $\mcF$ of $\mcZ(M)$,
    \begin{equation}\label{eq:ftd}
      r(\cup \mcF) = \sum_{\mcF'\subseteq \mcF\,:\,\mcF'\ne\emptyset}
      (-1)^{|\mcF'|+1}r(\cap \mcF'),
    \end{equation}
  \item for all ideals $\mcF\subseteq\mcZ(M)$,
    \begin{equation*}
      \sum_{Y\in\mcF}\alpha(Y) = \eta(\cup \mcF).
    \end{equation*}  
  \end{enumerate}
\end{thm}

We now consider how the results above extend to transversal matroids
of finite rank on infinite sets.  Although the ground set is infinite,
every multiset we consider is finite.  Thus, let $M$ be $M[\mcA]$
where $\mcA = (A_1,A_2,\ldots,A_r)$ is a set system on the infinite
set $E(M)$.  For each subset $X$ of $E(M)$, let $\phi(X)=\{k\,:\,X\cap
A_k\ne\emptyset\}$.  It is easy to see that if $F$ is a cyclic flat of
$M$, then $F = \{x\,:\,\phi(x)\subseteq \phi(F)\}$.  It follows that
$M$ has at most $\binom{r}{k}$ cyclic flats of rank $k$, so $\mcZ(M)$
is a finite lattice.  Whenever $M$ has finite rank and $\mcZ(M)$ is
finite, the definition of $\beta$ makes sense, as do the sums that
appear in the results above.  Reviewing the proofs shows that Theorems
\ref{thm:betami}, \ref{thm:ingthm2}, \ref{thm:ftm}, \ref{thm:ftmbeta},
\ref{thm:ingthm3}, and Corollary \ref{cor:tisft} hold in this setting,
where we add to the hypotheses of all but the last result the
requirements that $M$ has finite rank and $\mcZ(M)$ is finite.  Note
that in this setting, the assertion that matroids with fundamental
bases are transversal holds since the argument proving statement (2)
in Theorem \ref{thm:betami} shows that such matroids satisfy that
statement (with equality).  In contrast, Theorem \ref{thm:dualftm} was
obtained by duality, which does not apply within the class of matroids
of finite rank on infinite sets.  However, we have the following
result.

\begin{thm}\label{thm:infdualftm}
  A matroid $M$ of finite rank on an infinite set is a fundamental
  transversal matroid if and only if the lattice $\mcZ(M)$ is finite
  and equation (\ref{eq:ftd}) holds for all of its subsets
  (equivalently, ideals; equivalently, antichains).
\end{thm}

\begin{proof}
  First assume $M$ is a fundamental transversal matroid.  Let $X$ be a
  finite subset of $E(M)$ whose subsets include a fundamental basis, a
  cyclic spanning set for each cyclic flat, and a spanning set for
  each intersection of cyclic flats.  It follows that $M|X$ is a
  fundamental transversal matroid and the map $\psi:
  \mcZ(M|X)\rightarrow \mcZ(M)$ given by $\psi(Y)=\cl_M(Y)$ is a
  rank-preserving isomorphism.  Since $M|X$ is fundamental, the
  counterpart of equation (\ref{eq:ftd}) holds for $M|X$.  By using
  $\psi$, we can deduce equation (\ref{eq:ftd}) for $M$.

  To prove the converse, let $r=r(M)$ and let $X$ be a finite subset
  of $E(M)$ that contains a cyclic spanning set for each cyclic flat
  and a spanning set for each intersection of cyclic flats.  As above,
  the map $\psi$ given by $\psi(Y)=\cl_M(Y)$ is a rank-preserving
  isomorphism of $\mcZ(M|X)$ onto $\mcZ(M)$.  Using $\psi$, from the
  validity of equation (\ref{eq:ftd}) for $M$ we can deduce its
  counterpart for $M|X$, so $M|X$ is fundamental by Theorem
  \ref{thm:dualftm}.  Thus, some injection $\phi:\mcZ(M|X)\rightarrow
  2^{[r]}$ satisfies properties (1)--(3) of Theorem \ref{thm:ingthm3}.
  Define $\phi':\mcZ(M)\rightarrow 2^{[r]}$ by $\phi'(F) = \phi(F|X)$.
  It is immediate that $\phi'$ is an injection that satisfies
  properties (1)--(3) of Theorem \ref{thm:ingthm3}, so $M$ is
  fundamental.
\end{proof}

Brylawski's characterization of fundamental transversal
matroids~\cite[Proposition 4.2]{aff}, which we state next, follows
easily from Theorem \ref{thm:dualftm}.

\begin{thm}
  A matroid $M$ is a fundamental transversal matroid if and only if
  for all families $\mcF$ of intersections of cyclic flats,
  \begin{equation}
    r(\cup \mcF) \geq \sum_{\mcF'\subseteq \mcF\,:\,\mcF'\ne\emptyset}
    (-1)^{|\mcF'|+1}r(\cap \mcF'),\label{tom} 
  \end{equation}
  or, equivalently, equality holds in inequality (\ref{tom}).  The
  same statement holds for matroids of finite rank on infinite sets
  where, in the second part, we add that $\mcZ(M)$ is finite.
\end{thm}

\begin{proof} 
  An inclusion-exclusion argument like that in the proof of Theorem
  \ref{thm:betami} shows that equality holds in inequality (\ref{tom})
  for fundamental transversal matroids.  For the converse, note that
  if that equality always holds, then $M$ is fundamental by Theorem
  \ref{thm:dualftm}.  Thus, it suffices to show that having inequality
  (\ref{tom}) hold for all families $\mcF$ of intersections of cyclic
  flats of $M$ yields equality.  To prove this, we induct on $|\mcF|$.
  The base case, $|\mcF|=1$, is obvious.  Assume that $|\mcF|>1$ and
  that equality holds for all families of intersections of cyclic
  flats that have fewer sets than $\mcF$.  Fix $X\in\mcF$ and let
  $\mcF_{\widehat{X}} = \mcF-\{X\}$.  The set $\mcF' = \{X\}$
  contributes $r(X)$ to the sum.  The sets $\mcF'$ with
  $\mcF'\subseteq \mcF_{\widehat{X}}$ give terms that, by the
  induction hypothesis, together contribute $r(\cup
  \mcF_{\widehat{X}})$ to the sum.  The sets $\mcF'$ with
  $\{X\}\subset \mcF'$ contribute terms that are the negatives of the
  terms in the corresponding sum based on the family $\{F\cap
  X\,:\,F\in \mcF_{\widehat{X}}\}$; by the induction hypothesis, the
  sum of these terms is $-r(X\cap(\cup \mcF_{\widehat{X}}))$.  Thus,
  inequality (\ref{tom}) is equivalent to
  \begin{equation*}
    r(\cup \mcF)\geq r(X) + r(\cup \mcF_{\widehat{X}}) -
    r(X\cap(\cup \mcF_{\widehat{X}})).
  \end{equation*}  
  Semimodularity (the opposite inequality) gives equality.  This
  completes the induction.
\end{proof}

Finally, we apply our results to the free product, which was
introduced and studied by Crapo and Schmitt~\cite{fpm,uft}.  Given
matroids $M$ and $N$ on disjoint sets, their \emph{free product}
$M\frp N$ is the matroid on $E(M)\cup E(N)$ whose bases are the sets
$B\subseteq E(M)\cup E(N)$ with (i) $|B|=r(M)+r(N)$, (ii) $B\cap E(M)$
independent in $M$, and (iii) $B\cap E(N)$ spanning $N$
(see~\cite[Proposition 3.3]{uft}).  In general, $M\frp N\ne N\frp M$.
Relative to the weak order, the free product is the greatest matroid
$M'$ on $E(M)\cup E(N)$ with $M'\del E(N)=M$ and $M'/E(M)=N$.  Special
cases of the free product include the free extension of $M$ (set
$N=U_{0,1}$) and the free coextension of $N$ (set $M=U_{1,1}$).  The
dual of the free product is given by $(M\frp N)^*=N^*\frp M^*$.  The
following result is~\cite[Proposition 6.1]{uft}.

\begin{prop}\label{prop:cffp}
  A subset $F$ of $E(M)\cup E(N)$ other than $E(M)$ is in $\mcZ(M\frp
  N)$ if and only if either (i) $F\subset E(M)$ and $F\in \mcZ(M)$ or
  (ii) $E(M)\subset F$ and $F-E(M)\in\mcZ(N)$.  The set $E(M)$ is in
  $\mcZ(M\frp N)$ if and only if $E(M)\in \mcZ(M)$ and
  $\emptyset\in\mcZ(N)$.
\end{prop}

By giving a presentation of $M\frp N$ from presentations of $M$ and
$N$, Crapo and Schmitt \cite[Proposition 4.14]{uft} showed that free
products of transversal matroids are transversal.  The following
extension of their result can be proven using either ideas in
\cite{uft} or, as we sketch below, Theorems \ref{thm:betami} and
\ref{thm:ftm}.

\begin{thm}
  For matroids $M$ and $N$ on disjoint ground sets, their free product
  $M \frp N$ is transversal if and only if both $M$ and $N$ are.  The
  corresponding statements hold for fundamental transversal matroids,
  for cotransversal matroids, and for matroids that are both
  transversal and cotransversal.
\end{thm}

This result follows from Theorems \ref{thm:betami} and \ref{thm:ftm}
(using antichains) and three observations:
\begin{enumerate}
\item by Proposition \ref{prop:cffp}, any antichain in $\mcZ(M\frp N)$
  is either (i) an antichain in $\mcZ(M)$ or (ii) formed from an
  antichain in $\mcZ(N)$ by augmenting each set by $E(M)$;
\item if $X\subseteq E(M)$, then $r_{M\frp N}(X)=r_M(X)$ since $(M\frp
  N)\del E(N)=M$;
\item if $X\subseteq E(N)$, then $r_{M\frp N}(X\cup E(M))=r_N(X)+r(M)$
  since $(M\frp N)/E(M)$ is $N$.
\end{enumerate}

\end{document}